\documentclass{amsart}
\usepackage{amsmath,amssymb,etoolbox,mathrsfs,booktabs,float,enumerate,calc,graphicx}

\theoremstyle{plain}
\newtheorem{theorem}{Theorem}

\newtheorem*{proposition*}{Proposition}

\newtheorem{lemmma*}{Lemma}

\newtheorem{conjecture}[theorem]{Conjecture}

\theoremstyle{definition}

\newtheorem*{definition*}{Definition}

\newtheorem{remark*}{Remark}

\newcommand{\Irr}{\mathrm{Irr}}

\newcommand{\II}{{I\!I}}
\begin{document}
\title[Zeros in character tables of symmetric groups]{Large-scale Monte Carlo simulations for zeros in character tables of symmetric groups}
\author[A.\ R.\ Miller]{Alexander Rossi Miller}
\author[D.\ Scheinerman]{Danny Scheinerman}
\address{Center for Communications Research, Princeton, NJ}
\begin{abstract}
This is a brief report on some recent large-scale Monte Carlo simulations for
approximating the density of zeros in character tables of large symmetric groups.
Previous computations suggested that a large fraction of zeros
cannot be explained by classical vanishing results.
Our computations eclipse previous ones and suggest that the opposite is true. In fact, we find empirically that almost all of the zeros are of a single classical type.
\end{abstract}
\maketitle
\section{Introduction}
\subsection{A hundred years of computing}
For more than a century now, there has been much interest in
computing character values of $S_n$.
In 1900, Frobenius~\cite{Frobenius} discovered a formula that
allows one to compute by hand the full character table of $S_n$ when $n$ is small.
By 1934, Littlewood and Richardson~\cite{LR} had worked out the tables up to $n=9$.
A year later, Littlewood~\cite{Littlewood} finished the table for $n=10$ and Zia-ud-Din~\cite{Z1}
had the table for $n=11$. Two years after that, Zia-ud-Din~\cite{Z2}
worked out $n=12$ and $n=13$. Then in the late 1930's and early 1940's a new recursive method
for computing character values of $S_n$ came into play. It was introduced by Murnaghan and Nakayama~\cite{Murnaghan,Nakayama}, and it went on to become the preferred method for character calculations, starting with Kond\^{o}'s computation~\cite{Kondo} of the character table for $n=14$ in 1940. Already for $n=14$, the character table has over 18000 entries. 
About a decade later, in 1954, electronic computers entered the picture when, under the auspices of the Atomic Energy Commission at Los Alamos National Lab, 
Bivins--Metropolis--Stein--Wells~\cite{BMSW}
used the now-famous MANIAC I to compute the tables up to $n=16$ (over 53000 entries).
The MANIAC (short for Mathematical Analyzer Numerical Integrator and Automatic Computer) was the most important and powerful computer of its day. Around 5 years later, Stig Com\'et~\cite{Comet} used Sweden's first electronic computer, named BESK (the Binary Electronic Sequence Calculator), to compute many entries, but not all of them, in the character tables up to $n=37$. Com\'et went on to head the Swedish Board for Computing Machinery, the team that developed BESK. 
Nowadays, modern computer algebra software such as Magma and GAP can compute the complete tables up to around $n=37$ or $n=38$ ($\approx 6.7\times 10^8$ entries) on a typical good home computer, with the limitations being time and memory. But why keep computing?\pagebreak

\subsection{Statistical properties hiding in plain sight}
About ten years ago, the first author~\cite{Miller2014} started investigating statistical properties of character tables and discovered some remarkable results. 
Miller~\cite{Miller2014} considered random character values $\chi(g)$ of a finite group $G$ and  asked,
\[\text{\emph{What is the chance that ${\chi(g)=0}$?}}\]

Here there are two natural ways of choosing a character value at random. The first way is to choose the character $\chi$ uniformly at random from the set of irreducible characters ${\rm Irr}(G)$, choose the group element $g$ uniformly at random from the group $G$, and then evaluate $\chi(g)$. The main result of~\cite{Miller2014} is then that for $G=S_n$, $\chi(g)=0$ with probability $\to 1$ as $n\to\infty$. That is,
\[\frac{|\{(\chi,g)\in {\rm Irr}(S_n)\times S_n : \chi(g)=0\}|}{|\Irr(S_n)\times S_n|}\to 1\text{ as $n\to\infty$}.\]

The second way to choose a character value $\chi(g)$ is to simply choose an entry uniformly at random from the character table of $G$. That is, choose $\chi$ uniformly at random from ${\rm Irr}(G)$, choose a class $g^G$ uniformly at random from the set of all conjugacy classes ${\rm Cl}(G)$, and then evaluate $\chi(g)$. Now the chance that {$\chi(g)$ equals~$0$} is the fraction of the character table covered by zeros, i.e.\ the fraction
\begin{equation}\label{zeros density}
  z(G)=\frac{|\{(\chi,g^G)\in {\rm Irr}(G)\times {\rm Cl}(G) : \chi(g)=0\}|}{|{\rm Irr}(G)\times {\rm Cl}(G)|}.
\end{equation}

\begin{figure}[H]
    \includegraphics[scale=1.1]{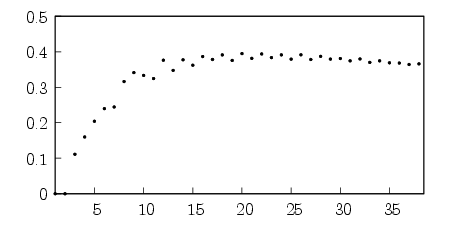}
    \caption{$z(S_n)$ for $1\leq n\leq 38$.}\label{small zeros z}
  \end{figure}

  The numerator of $z(S_n)$ counts the total number of zeros in
  the character table of $S_n$, a number which was first computed for
  various values of $n$ by John McKay and communicated to N. J. A. Sloane
  around 1991. The sequence was extended in 2012 to $n=35$ by Eric M.\ Schmidt.
  Miller provided the numbers up to $n=38$ in~\cite{Miller2019}. 
  Dividing these numbers by the appropriate denominator in \eqref{zeros density},
  we obtain the plot for $z(S_n)$ with $1\leq n\leq 38$ shown in Figure~\ref{small zeros z}.

  From the available data up to $n=38$, the limiting behavior of $z(S_n)$
  is anyone's guess. Does it tend to $0$? Can it be bounded away from $0$?
  Does the limit even exist?  One can guess,
  but the fact of the matter is that we simply do not have enough data.

  \subsection{New computations and conjectures}\label{S New Comp}  The goal of the present paper is to shed some light on this problem
  by carrying out large-scale Monte Carlo simulations.
  Our findings provide valuable  new insights about
  the bulk of zeros in very large character tables of symmetric groups.
  
Our Monte Carlo simulations could not be done with available packages such as GAP and Magma. To demonstrate, we generated 1000 pairs $(\chi,g^G)$ uniformly at random from ${\Irr}(G)\times {\rm Cl}(G)$ for  $G=S_{47}$, stored them in a file, and then proceeded to evaluate $\chi(g)$ using Magma and GAP on a 3GHz machine. Magma took 17863 seconds ($\approx 5$ hours) and GAP took 171173 seconds ($\approx 48$ hours). On average, the run-times get exponentially worse as $n$ increases.
But a more serious problem is that soon after $n=47$, both Magma and GAP either abort with exorbitant memory requests (Magma aborted one trial for $n=48$ with a reported failed memory request of 17592186042368.1MB) or, after only a handful of trials, the algebra software gets stuck on a single computation with no end in sight. 

Not only were we aiming to go beyond $n=47$, but we were also aiming for many more than a thousand samples for each $n$. To solve this problem, we wrote
a small suite in the Julia language and made use of various performance-enhancing techniques. We briefly comment on these techniques at the end of the paper. 
Our program took around $0.314$ seconds to evaluate the same 1000 different values $\chi(g)$ for $S_{47}$ that took GAP and Magma several hours, all else being equal. So it was a little more than 540000 times faster than GAP and more than 56000 times faster than Magma for this random batch of entries in the character table of $S_{47}$. This put our target computations within reach of just a few days on a modern multi-core chip with ample random-access memory.

With our program we were able to go well beyond $n=47$ and draw many more than a thousand samples for each $n$.
These new computations are far out of reach of modern computer algebra packages, and the resulting data provides key new insights about the zeros in large character tables of symmetric groups.

\section{New data and conjectures}
\subsection{Notation and conventions}
By a partition $\lambda$ of a positive integer $n$, we mean
a sequence $\lambda=(\lambda_1,\lambda_2,\ldots,\lambda_{\ell(\lambda)})$ with parts $\lambda_k$
that are positive integers satisfying 
$\lambda_1\geq \lambda_2\geq \ldots\geq \lambda_{\ell(\lambda)}$ and $\sum_k \lambda_k=n$.
To each pair of positive integers $(i,j)$ satisfying
$1\leq i\leq \ell(\lambda)$ and $1\leq j\leq \lambda_i$, there is the  hook length \[h_\lambda(i,j)=\lambda_i-j+|\{s\in\{i,i+1,\ldots,\ell(\lambda)\} : \lambda_s\geq j\}|.\] The hook lengths $h_\lambda(i,j)$, for all admissible pairs $(i,j)$, constitute the hook lengths of $\lambda$. If none of the hook lengths of $\lambda$ are divisible by a given positive integer $t$, then $\lambda$ is called a $t$-core partition.

For a partition $\lambda$ of a positive integer $n$, we  write $\chi_\lambda$ for the irreducible character of $S_n$ associated with $\lambda$ in the usual way. By $\chi_\lambda(\mu)$, we  mean the value $\chi_\lambda(g)$ for any permutation $g\in S_n$ whose cycle type is $\mu$ in the sense that the periods of the various cycles for $g$ make up the parts of $\mu$. The character table of $S_n$ is then the table $[\chi_\lambda(\mu)]_{\lambda,\mu}$ with $\lambda$ and $\mu$ running over all partitions of $n$. We denote by $P_n$ the set of all partitions of $n$, and we denote by $p_n$ the number of partitions of $n$, so there are
$p_n^2=|{\rm Irr}(S_n)\times {\rm Cl}(S_n)|$
entries in the character table, and the density
of zeros in the table is
\[z(S_n)=\frac{|\{(\lambda,\mu)\in P_n^2 : \chi_\lambda(\mu)=0\}|}{p_n^2}.\]

\subsubsection{Types of zeros}
Let
\[\mathcal Z(S_n)=\{(\lambda,\mu)\in P_n^2 : \chi_\lambda(\mu)=0\}.\]
For the purposes of this paper, it will be convenient to
simply refer to a pair
$(\lambda,\mu)\in \mathcal Z(S_n)$ as a zero, instead of ``the location of a zero in the character table''.
It can be difficult to determine if a given pair $(\lambda,\mu)$ belongs to $\mathcal Z(S_n)$, but there are two well-known special cases that are  easier to handle (see e.g.~\cite{JamesKerber, Miller2014,Miller2019,PS1,McSpiritOno,PS2}).
\begin{enumerate}[(I)]
\item If $\lambda$ is a $\mu_1$-core, then $(\lambda,\mu)\in\mathcal Z(S_n)$ and we call $(\lambda,\mu)$ a zero of \emph{type I}.
\item 
If $\lambda$ is a $\mu_k$-core for some part $\mu_k$ of $\mu$,
then $(\lambda,\mu)\in \mathcal Z(S_n)$ and we call $(\lambda,\mu)$ a zero of \emph{type II}.
\end{enumerate}
We let $\mathcal Z_I(S_n)$ (resp.\ $\mathcal Z_\II(S_n)$) denote the set of all type I
(resp.\ type II) zeros in $P_n^2$, so
\[\mathcal Z_I(S_n)\subseteq \mathcal Z_\II(S_n)\subseteq \mathcal Z(S_n).\]
Let $z_I(S_n)$ (resp.\ $z_\II(S_n)$) denote the fraction of pairs $(\lambda,\mu)\in P_n^2$ that belong to $\mathcal Z_I(S_n)$ (resp.\ $\mathcal Z_\II(S_n)$), so
\[z_I(S_n)=\frac{|\mathcal Z_I(S_n)|}{p_n^2},\quad
  z_\II(S_n)=\frac{|\mathcal Z_\II(S_n)|}{p_n^2},\quad
  z(S_n)=\frac{|\mathcal Z(S_n)|}{p_n^2},\]
and
\[z_I(S_n)\leq z_\II(S_n)\leq z(S_n).\]

\begin{figure}[H]
  \includegraphics{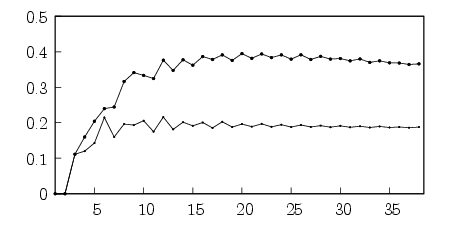}
  \caption{$z(S_n)$ and the lower bound $z_I(S_n)$ for $1\leq n\leq 38$.}\label{compare z and zii up to 38}
\end{figure}

Using the exact values available up to $n=38$, we have the plot shown in Figure~\ref{compare z and zii up to 38} for both $z_I(S_n)$ and $z(S_n)$. The ratios $z_I(S_n)/z(S_n)$ are given in Table~\ref{table}.

This comparison of $z_I(S_n)$ and $z(S_n)$ up to $n=38$ seems to suggest that many of the zeros in the character table, roughly half of them, are not of the first type, and therefore must be explained by some new, more complicated vanishing result.
\begin{center}
  \begin{table}[H]
\begin{tabular}{lcclcclc}
  \toprule
  $n$ & $\frac{z_I(S_n)}{z(S_n)}$ &&
                             $n$ & $\frac{z_I(S_n)}{z(S_n)}$ &&
  $n$ & $\frac{z_I(S_n)}{z(S_n)}$\\
  \cmidrule{1-2}\cmidrule{4-5}\cmidrule{7-8}
  $3$  & 1.000 & & 15 & 0.529 & & 27 & 0.497 \\
  $4$  & 0.750 & & 16 & 0.519 & & 28 & 0.496 \\
  $5$  & 0.700 & & 17 & 0.491 & & 29 & 0.494 \\
  $6$  & 0.897 & & 18 & 0.518 & & 30 & 0.502 \\
  $7$  & 0.655 & & 19 & 0.500 & & 31 & 0.500 \\
  $8$  & 0.621 & & 20 & 0.497 & & 32 & 0.500 \\
  $9$  & 0.567 & & 21 & 0.495 & & 33 & 0.504 \\
  $10$ & 0.617 & & 22 & 0.500 & & 34 & 0.506 \\
  $11$ & 0.538 & & 23 & 0.491 & & 35 & 0.506 \\
  $12$ & 0.574 & & 24 & 0.497 & & 36 & 0.511 \\
  $13$ & 0.522 & & 25 & 0.495 & & 37 & 0.511 \\
  $14$ & 0.534 & & 26 & 0.495 & & 38 & 0.513 \\
  \bottomrule\\
\end{tabular}\\ 
    \caption{\label{table} The fraction $z_I(S_n)/z(S_n)$
    for $3\leq n\leq 38$, rounded to the number of decimal places shown.}
\end{table}
\end{center}
\subsection{Monte Carlo simulations and conjectures}

\subsubsection{One million samples for each $n$ up to 150}
\begin{figure}[H]
  \includegraphics[width=\textwidth]{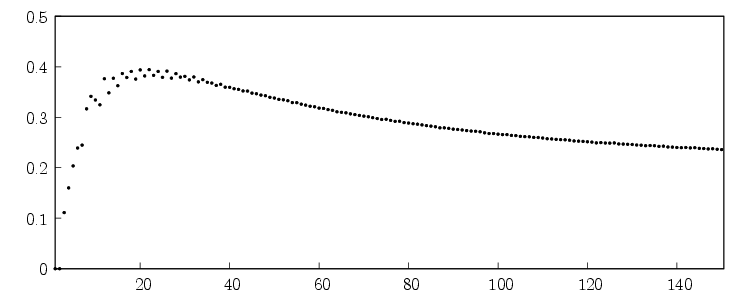}
  \caption{Approximation of $z(S_n)$ for $1\leq n\leq 150$
  based on $10^6$ simulations for each $n$.}\label{MC 1 million samples}
\end{figure}

For our first Monte Carlo experiment, we drew 
$10^6$ pairs $(\lambda,\mu)$ uniformly at random from $P_n^2$ for each $n$ up to $150$.
Figure~\ref{MC 1 million samples} shows the plot for the fraction of
these samples $(\lambda,\mu)$ such that $\chi_\lambda(\mu)=0$.
This offers the first convincing empirical evidence that
$\lim_{n\to\infty}z(S_n)$ indeed exists.

\subsubsection{Comparing  $z_{I}(S_n)$, $z_{\II}(S_n)$, and $z(S_n)$}
Perhaps the most important outcome of our Monte Carlo simulations is
that, in contrast to the old picture up to $n=38$, for large $n$ we find that 
 almost all zeros are at least of the second type.

\begin{figure}[H]
  \includegraphics[width=\textwidth]{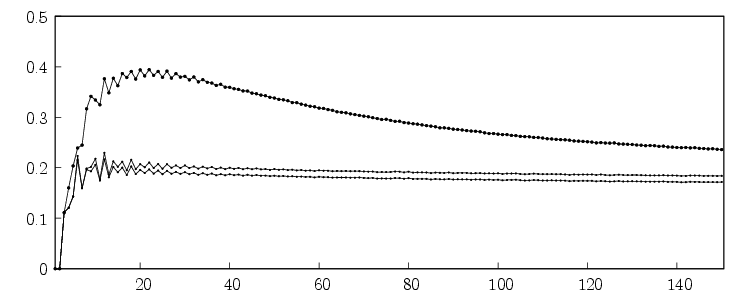}
  \caption{Approximations of $z_{I}(S_n),z_{\II}(S_n),z(S_n)$~for~${1\leq n\leq 150}$  based on $10^6$ random samples $(\lambda,\mu)\in P_n^2$ for each $n$.}\label{MC all types 1 million}
\end{figure}

Figure~\ref{MC all types 1 million} shows, for $1\leq n\leq 150$,
the following three fractions of $10^6$ samples $(\lambda,\mu)\in P_n^2$, in decreasing order:
   the fraction of samples with $(\lambda,\mu)\in \mathcal Z(S_n)$,
   the fraction with $(\lambda,\mu)\in \mathcal Z_\II(S_n)$, and
    the fraction with $(\lambda,\mu)\in\mathcal Z_I(S_n)$.

Table~\ref{table approx 10^6} records the approximations
of $z_I(S_n),z_\II(S_n),z(S_n)$ for $n=10,20,\ldots,150$.
Here and in subsequent tables, our approximations for $z_I(S_n),z_\II(S_n),z(S_n)$ 
are rounded to the number of decimal places shown.

\begin{center}%
  \begin{table}[H]%
\begin{tabular}{lcccclccc}
\toprule
$n$ & $z_I(S_n)$ & $z_\II(S_n)$ & $z(S_n)$ & & $n$ & $z_I(S_n)$ & $z_\II(S_n)$ & $z(S_n)$\\
\cmidrule{1-4}\cmidrule{6-9}%
5  & 0.143 & 0.143 & 0.204 & & 80  & 0.179 & 0.192 & 0.289 \\
10 & 0.206 & 0.218 & 0.334 & & 85  & 0.177 & 0.189 & 0.282 \\
15 & 0.191 & 0.202 & 0.362 & & 90  & 0.177 & 0.189 & 0.276 \\
20 & 0.196 & 0.207 & 0.394 & & 95  & 0.177 & 0.189 & 0.272 \\
25 & 0.187 & 0.198 & 0.379 & & 100 & 0.176 & 0.189 & 0.267 \\
30 & 0.191 & 0.205 & 0.381 & & 105 & 0.175 & 0.188 & 0.262 \\
35 & 0.186 & 0.199 & 0.369 & & 110 & 0.175 & 0.187 & 0.259 \\
40 & 0.187 & 0.200 & 0.359 & & 115 & 0.174 & 0.187 & 0.255 \\
45 & 0.184 & 0.197 & 0.348 & & 120 & 0.174 & 0.186 & 0.252 \\
50 & 0.184 & 0.197 & 0.338 & & 125 & 0.173 & 0.185 & 0.249 \\
55 & 0.182 & 0.195 & 0.329 & & 130 & 0.173 & 0.186 & 0.246 \\
60 & 0.182 & 0.195 & 0.318 & & 135 & 0.173 & 0.185 & 0.244 \\
65 & 0.181 & 0.193 & 0.310 & & 140 & 0.172 & 0.184 & 0.240 \\
70 & 0.180 & 0.193 & 0.302 & & 145 & 0.172 & 0.184 & 0.238 \\
75 & 0.179 & 0.191 & 0.296 & & 150 & 0.172 & 0.184 & 0.236 \\
\bottomrule\\
\end{tabular}\\ 
\caption{
  Approximations of $z_I(S_n),z_\II(S_n),z(S_n)$ based
  on $10^6$  samples $(\lambda,\mu)\in P_n^2$ for
  each $n$.\label{table approx 10^6}
  }
\end{table}
\end{center}

\begin{figure}[H]
  \includegraphics[width=\textwidth]{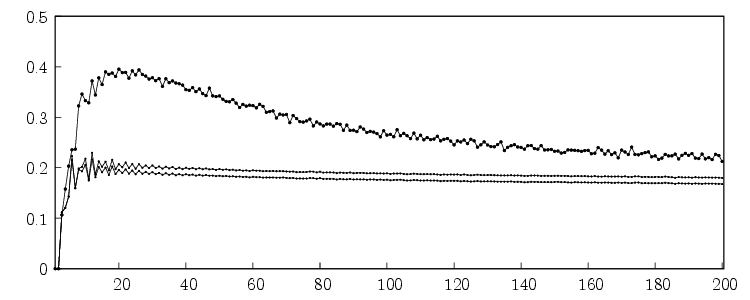}
  \caption{Approximations of $z_{I}(S_n),z_{\II}(S_n),z(S_n)$~for~${1\leq n\leq 200}$  based on $10^4$ samples for
    $z(S_n)$ and $10^6$ samples for $z_I(S_n),z_\II(S_n)$.}\label{MC all types up to 200}
\end{figure}

We also ran the same Monte Carlo experiment at a lower resolution
for $z(S_n)$ using $10000$ samples (instead of $10^6$) for each $n$ up to $200$.
We kept the resolution for $z_I(S_n)$ and $z_\II(S_n)$ at $10^6$ samples.
The results are shown in Figure~\ref{MC all types up to 200}. (See also Table~\ref{t3}.)

\begin{center}%
  \begin{table}[H]%
\begin{tabular}{lccc}
  \toprule
  $n$& $ z_I(S_n)$ & $ z_\II(S_n)$ & $ z(S_n)$ \\
  \cmidrule{1-4}%
150 & 0.172 & 0.184 & 0.233 \\
155 & 0.171 & 0.183 & 0.235 \\
160 & 0.171 & 0.183 & 0.235 \\
165 & 0.171 & 0.183 & 0.227 \\
170 & 0.170 & 0.182 & 0.235 \\
175 & 0.171 & 0.183 & 0.226 \\
180 & 0.169 & 0.182 & 0.224 \\
185 & 0.169 & 0.181 & 0.224 \\
190 & 0.169 & 0.181 & 0.225 \\
195 & 0.169 & 0.181 & 0.218 \\
200 & 0.168 & 0.180 & 0.212 \\
  \bottomrule\\
\end{tabular}\\ 
\caption{\label{t3}
  Approximations of $z_I(S_n),z_\II(S_n),z(S_n)$ based
  on $10^4$ samples for $z(S_n)$ and $10^6$ samples for $z_I(S_n)$ and
  $z_\II(S_n)$.
    }
\end{table}
\end{center}

\subsubsection{The difference $z_\II(S_n)-z_I(S_n)$}
\begin{figure}[H]%
\includegraphics[width=\textwidth]{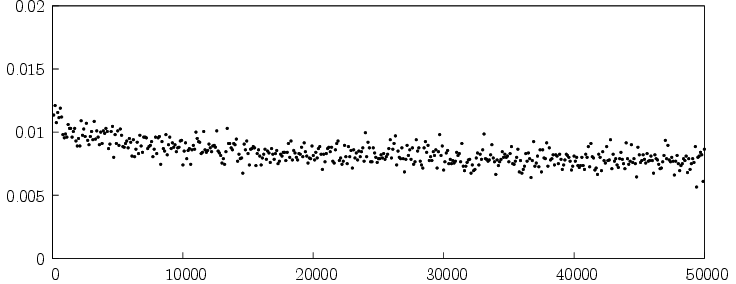}%
\caption{Approximation of the difference $z_\II(S_n)-z_I(S_n)$ for
    $n=100,200,300,\ldots,50000$ 
    based on  $20000$ samples $(\lambda,\mu)\in P_n^2$ for each $n$.}\label{MC difference trivial}
\end{figure}

To understand better the difference between $z_\II(S_n)$ and $z_I(S_n)$,
we drew $20000$ samples for $n=100,200,300,\ldots, 50000$
and looked at the fraction of the samples that were of type II but not type I. 
The results are shown in Figure~\ref{MC difference trivial}. See also Table~\ref{zI and zII approx}. 
We suspect that the difference $z_\II(S_n)-z_I(S_n)$ tends to zero and is in fact $o(z_I(S_n))$.

\subsubsection{Approximating $z_I(S_n)$}\label{z_I section}
\begin{figure}[H]%
\includegraphics[width=\textwidth]{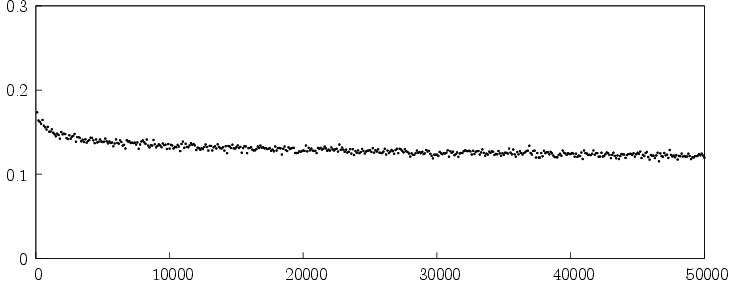}%
\caption{Approximation of $z_{I}(S_n)$ for
    ${n=100,200,\ldots,50000}$ 
    based on $20000$ samples $(\lambda,\mu)\in P_n$ for each $n$.}\label{figure type 1 50k}
\end{figure}

Lastly, for $n=100,200,300,\ldots,50000$,
we 
drew $20000$ samples $(\lambda,\mu)\in P_n^2$
and looked at the fraction that were of type I.
The results are shown in
Figure~\ref{figure type 1 50k}.
After circulating this paper, Peluse and Soundararajan communicated 
to the authors an unpublished result of theirs that
\[z_I(S_n)\sim 2/\log n,\]
and independently 
Ryan Eberhart used symbolic regression and 
the values of $z_I(S_n)$ that we had computed up to around $n=5000$
to experimentally discover the same result up to a constant factor.

\begin{center}%
  \begin{table}[H]%
\begin{tabular}{lccclcc}
  \toprule
  $n$& $z_I(S_n)$ & $z_\II(S_n)$
  && $n$& $z_I(S_n)$ & $z_\II(S_n)$\\
  \cmidrule{1-3}\cmidrule{5-7}
1000 & 0.150 & 0.160 && 26000 & 0.127 & 0.135 \\
2000 & 0.147 & 0.157 && 27000 & 0.130 & 0.137 \\
3000 & 0.139 & 0.149 && 28000 & 0.121 & 0.129 \\
4000 & 0.140 & 0.150 && 29000 & 0.124 & 0.132 \\
5000 & 0.138 & 0.148 && 30000 & 0.123 & 0.130 \\
6000 & 0.141 & 0.151 && 31000 & 0.127 & 0.136 \\
7000 & 0.139 & 0.149 && 32000 & 0.127 & 0.135 \\
8000 & 0.140 & 0.148 && 33000 & 0.127 & 0.136 \\
9000 & 0.132 & 0.141 && 34000 & 0.125 & 0.132 \\
10000 & 0.130 & 0.138 && 35000 & 0.123 & 0.131 \\
11000 & 0.139 & 0.149 && 36000 & 0.126 & 0.133 \\
12000 & 0.129 & 0.138 && 37000 & 0.125 & 0.133 \\
13000 & 0.134 & 0.141 && 38000 & 0.128 & 0.135 \\
14000 & 0.132 & 0.141 && 39000 & 0.120 & 0.128 \\
15000 & 0.132 & 0.140 && 40000 & 0.124 & 0.132 \\
16000 & 0.131 & 0.139 && 41000 & 0.128 & 0.136 \\
17000 & 0.132 & 0.140 && 42000 & 0.125 & 0.133 \\
18000 & 0.128 & 0.137 && 43000 & 0.119 & 0.127 \\
19000 & 0.130 & 0.138 && 44000 & 0.123 & 0.132 \\
20000 & 0.128 & 0.136 && 45000 & 0.126 & 0.134 \\
21000 & 0.125 & 0.133 && 46000 & 0.123 & 0.131 \\
22000 & 0.128 & 0.136 && 47000 & 0.125 & 0.135 \\
23000 & 0.129 & 0.136 && 48000 & 0.117 & 0.125 \\
24000 & 0.125 & 0.135 && 49000 & 0.118 & 0.125 \\
25000 & 0.129 & 0.136 && 50000 & 0.120 & 0.128 \\
  \bottomrule\\
\end{tabular}\\ 
\caption{
  Approximations of $z_I(S_n)$ and $z_\II(S_n)$ based
  on $20000$ random samples $(\lambda,\mu)\in P_n^2$ for
  each $n$.
  \label{zI and zII approx}
  }
\end{table}
\end{center}

\vspace{-.25in}
\subsubsection{Conjectures}
\begin{conjecture}\label{conj equality}
 $z(S_n)=O(\frac{1}{\log{n}})$ as $n\to\infty$.
\end{conjecture}

As remarked in \S\ref{z_I section}, 
$z_I(S_n)\sim\frac{2}{\log n}$,
so Conjecture~\ref{conj equality} is equivalent to the statement 
that $z(S_n)=O(z_I(S_n))$. We suspect that in fact
$z(S_n)$ is asymptotic to  $z_I(S_n)$, i.e.\ almost every zero in the character
table of $S_n$ is of type I as $n\to\infty$.

\begin{conjecture}\label{conj asym}
  $z(S_n)\sim z_I(S_n)\sim\frac{2}{\log n}$ as $n\to\infty$.
\end{conjecture}

\begin{conjecture}\label{M conj}
  $z_I(S_n)\geq z_I(S_{n+1})$ for $n\geq 82$.
\end{conjecture}
We have checked Conjecture~\ref{M conj} up to $n=5000$. The number of type I zeros for $S_{5000}$ is already $\approx 4\times 10^{147}$.  The exact number is\\
$|\mathcal Z_I(S_{5000})|=\,$401646560041542513478712635316593642393796744588210221412396\\202065981953145283377163249198113710031345332410560064801718418287227587\\3569544245365379.

\section{Comments on computational ingredients}
We end by briefly describing a couple of the main ingredients that made
our computations possible.
Although the methods and procedures are not new, they result in a tremendous speedup over existing software packages, as seen in~\S\ref{S New Comp}.

Just as with GAP and Magma, our method for computing a
particular character value $\chi_\lambda(\mu)$ is based on the Murnaghan--Nakayama rule,
which says that
\begin{equation}\label{MN}
  \chi_\lambda(\mu)=\sum (-1)^{{\rm ht}(\rho)}\chi_{\lambda\backslash \rho}(\nu)
\end{equation}
where $\nu=(\mu_2,\mu_3,\ldots)$, the sum is over all so-called \emph{rim hooks}
$\rho$ in $\lambda$ of size $\mu_1$, and ${\rm ht}(\rho)$ is one less than the number of rows occupied by $\rho$. See \cite{JamesKerber} for details.
By repeatedly expanding summands with the same rule, this gives a recursive procedure for evaluating $\chi_\lambda(\mu)$.

The recursive procedure for evaluating $\chi_\lambda(\mu)$ was carried out using a loop and, at each stage of the process, a dictionary structure was used to store the unique shapes $\lambda\backslash \rho$ that occur as indices, along with the  corresponding coefficients.  During this process, when we reach terms of the shape $\chi_{\xi}(\omega)$ with $\omega=(1,1,\ldots,1)$, we stop and evaluate the terms using the usual hook-length formula \cite{JamesKerber}.

\begin{figure}[H]
  \includegraphics[scale=1.2]{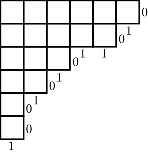}
  \caption{}\label{653211}
\end{figure}

Our primary speedup comes from encoding $\lambda$ as an unsigned integer and performing fast bit manipulations for going from a given partition $\lambda$ to the various partitions $\lambda\backslash \rho$ and associated signs that occur on the right-hand side of \eqref{MN}. 
The encoding of $\lambda$ as an integer is based on the usual binary encoding.
Viewing $\lambda$ as a Young diagram and walking along the outer boundary edges from the lower left to the upper right, we obtain a binary representation by recording a 1 for each horizontal edge and a 0 for each vertical edge. For example, the partition $(6,5,3,2,1,1)$
is encoded by the integer \texttt{0b100101011010}, which is the usual binary shorthand for the integer $2^{11}+2^8+2^6+2^4+2^3+2$. See Figure~\ref{653211}.
In this language of unsigned integers, the removal of a rim hook $\rho$ in $\lambda$ of size $\mu_1$ translates into flipping two bits: a $1$ and a $0$ lying $\mu_1$ steps to the right of the~$1$. The sign $(-1)^{{\rm ht}(\rho)}$ becomes $-1$ raised to the number of zeros that lie strictly between the two bits being flipped. For example, if $\lambda=(6,5,3,2,1,1)$ as before, so that we are working with $\lambda$ as the unsigned integer 
\texttt{0b100101011010}, and if $\mu_1=3$, then the resulting shapes $\lambda\backslash \rho$ that show up in \eqref{MN} are represented by \texttt{0b100001111010} and \texttt{0b100101010011}, both with a sign of $-1$. We perform these operations with fast bit manipulations.
\pagebreak

\subsubsection*{Acknowledgments}
It is a pleasure to thank James Barker, Reinier Br\"{o}ker, and Mike Mossinghoff
for their assistance with computational resources that we used to carry out our computations.
We also thank Leighton P.\ Barnes, Tim Chow, Ryan Eberhart, Charles Tomlinson, and Jeff VanderKam for discussions at various stages of the project, and Sarah Peluse and Kannan Soundararajan for sharing with us the unpublished result of theirs in \S\ref{z_I section}.

\end{document}